\documentclass[11pt,final]{article}
\usepackage{amsmath,latexsym,amssymb,amsmath,
amscd,amsthm,amsxtra}

\newtheorem{theorem}{Theorem}[section]

\def\RR{\mathbb{R}}

\def\SS{\mathbb{S}}
\def\HH{\mathbb{H}}
\def\tr{\mathrm{tr}}

\def\O{\Omega}

\def\De{\Delta}
\def\n{\nabla}
\def\ode{\overline{\Delta}}
\def\on{\overline{\nabla}}

\def\<{\langle}
\def\>{\rangle}

\begin{document}

\title{A generalization of Reilly's formula and its applications to a new Heintze-Karcher type inequality}
\author{Guohuan Qiu \footnote{Department of Mathematics, University of Science and Technology of China,  Hefei , P. R. China \quad Email: guohuan@mail.ustc.edu.cn.}    \and Chao Xia \footnote{Max-Planck-Institut f\"ur Mathematik in den Naturwissenschaften, Inselstr. 22, D-04103, Leipzig, Germany \quad Email: chao.xia@mis.mpg.de.}}

\date{}
\maketitle

\medskip

\abstract{In this paper, we prove a generalization of Reilly's formula in \cite{Reilly}. We apply such general Reilly's formula to give alternative proofs of the Alexandrov's Theorem and the Heintze-Karcher  inequality in the hemisphere and in the hyperbolic space. Moreover, we use the general Reilly's formula to prove a new Heintze-Karcher inequality for Riemannian manifolds with boundary and sectional curvature bounded below.}

\


\

\section{Introduction}

In a celebrated paper \cite{Reilly}, Reilly proved an integral formula for compact Riemannian manifolds with smooth boundary. To be precise, let us first give our notations. Throughout this paper, Let $(\Omega^{n}, g)$ be an $n$-dimensional compact  Riemannian manifold with smooth boundary $M$. We denote by $\on$, $\ode$ and $\on^2$ the gradient, the Laplacian and the Hessian on $\O$ respectively, while by $\n$ and $\De$ the gradient and the Laplacian on $M$ respectively.  Let $\nu$ be the unit outward normal of $M$. We denote by $h(X,Y)=g(\on_X \nu, Y)$ and $H=\frac{1}{n-1}\tr_g h$ the second fundamental form and the (normalized) mean curvature (with respect to $\nu$) of $M$ respectively. Let $d\O$ and $dA$ be the canonical measure of $\O$ and $M$ respectively. Let $\rm{Sect}$ and $\rm{Ric}$ be the sectional curvature and the Ricci curvature tensor of $\O$ respectively.

 Given a smooth function $f$   on $\O$,  we denote $z=f|_M$ and $u=\on_\nu f$. Reilly's formula \cite{Reilly} states that
\begin{eqnarray}\label{Reilly}
&&\int_\O \{(\ode f)^2-|\on^2 f|^2-Ric(\on f, \on f)\}d\O\nonumber\\&&=\int_M \{2u\De z+(n-1)Hu^2+h(\n z, \n z)\}dA.
\end{eqnarray}

Reilly's formula \eqref{Reilly} has numerous applications. For example, in \cite{Reilly} Reilly himself applied it to prove a Lichnerowicz type sharp lower bound for the first eigenvalue of the Laplacian on manifolds with boundary and reprove Alexandrov's rigidity theorem for  embedded hypersurfaces with constant mean curvature in $\mathbb{R}^n$. Other applications can be found for instance in \cite{CW, PRS, Ros}.

In \cite{Ros} Ros used Reilly's formula  to prove the following integral inequality, which was applied to show Alexandrov's rigidity theorem for high order mean curvatures. 
\vspace{3mm}

\noindent{\bf Theorem A.} (Ros \cite{Ros}) {\it Let $(\Omega^{n}, g)$ be a compact $n$-dimensional Riemannian manifold with smooth boundary $M$ and non-negative Ricci curvature. Let $H$ be the  mean curvature of $M$. If $H$ is positive everywhere, then
\begin{eqnarray}\label{ros}
\int_M \frac{1}{H}dA\geq n\rm{Vol}(\Omega).
\end{eqnarray}
The equality in \eqref{ros} holds if and only if $\O$ is isometric to an Euclidean ball.
}

\vspace{3mm}

For $\O\subset \RR^n$, inequality \eqref{ros} is essentially contained in the paper of Heintze and Karcher \cite{HK}. Ros' proof of Theorem A  based on the Reilly's formula \eqref{Reilly} and a suitable Dirichlet boundary value problem. Later, Montiel and Ros \cite{MR} gave an alternative proof of Theorem A in the case $\O\subset \mathbb{R}^n$ based on the ideas of Heintze and Karcher \cite{HK}, so that they had an alternative proof of Alexandrov's theorem in $\mathbb{R}^n$. Using the same idea of \cite{HK}, they also showed Alexandrov type theorem for constant higher order mean curvature embedded hypersurfaces in the hemi-sphere $\mathbb{S}_+^n$ or the hyperbolic space $\mathbb{H}^n$. However, they could not prove a similar inequality as \eqref{ros} in $\mathbb{S}_+^n$ or  $\mathbb{H}^n$. 

Quite recently, in order to study the Alexandrov type rigidity problem in general relativity,  Brendle \cite{Br} gave a version of Theorem A in $\mathbb{S}_+^n$ and $\mathbb{H}^n$, more generally, in a large class of warped product spaces, including the Schwarzschild manifold.

For $\mathbb{S}_+^n$ and $\mathbb{H}^n$, Brendle's result states as follows.
\vspace{3mm}

\noindent{\bf Theorem B.} (Brendle \cite{Br}) {\it Let $\Omega^{n}\subset \mathbb{H}^n$ ( $\SS_+^n$ resp.) be a compact $n$-dimensional domain with smooth boundary $M$. Let $H$ be the normalized mean curvature of $M$. Let $V(x)=\cosh \rm{dist}$$_{\HH^n}(x,0)$ ($\cos \rm{dist}$$_{\mathbb{S}^n}(x,0)$ resp.). If $H$ is positive everywhere, then
\begin{eqnarray}\label{brendle}
\int_M \frac{V}{H}dA\geq n\int_\Omega Vd\O.
\end{eqnarray}
The equality in \eqref{brendle} holds if and only if $\O$ is isometric to a geodesic ball.
}
\vspace{3mm}

As mentioned before, Brendle proved \eqref{brendle} for more general warped product spaces. Recently, his inequality has many interesting applications in general relativity, see for instance \cite{BHW, BW, dG}.

Brendle's method is quite different from Ros'. He used a geometric flow, along which the quantity $\int_M \frac{V}{H}dA$ is monotone non-increasing, to prove Theorem B. It is a natural problem to ask whether there is a Reilly-Ros type proof for Theorem B. This is the motivation of this paper.

\

In this paper,  we first prove the following general Reilly's formula.
\begin{theorem} Let $V: \overline{\O}\to \RR$ be a given a.e. twice differentiable function. Given a smooth function $f$  on $\O$,  we denote $z=f|_M$ and $u=\on_\nu f$. Let $K\in\RR$. Then we have the following identity:
\begin{eqnarray}\label{xeq0}
&&\int_\O V\left((\ode f+Knf)^2-|\on^2 f+Kfg|^2\right)d\O\nonumber\\&=&\int_M V\left(2u\De z+(n-1)Hu^2+h(\n z, \n z)+(2n-2)Kuz\right)dA\nonumber\\&&+\int_M \on_\nu V\left(|\n z|^2-(n-1)Kz^2\right) dA\nonumber\\&&+\int_\O \left(\on^2 V-\ode Vg-(2n-2)KVg+V Ric\right)(\on f, \on f)d\O\nonumber\\&&+(n-1)\int_\O (K\ode V+nK^2 V)f^2 d\O.
\end{eqnarray}
\end{theorem}
When $V\equiv 1$ and $K=0$, \eqref{xeq0} reduces to Reilly's formula \eqref{Reilly}. We are interested in some other choices of $V$ in this paper, particularly, $V(x)=\cosh r(x)$ or $\cos r(x),$ where $r(x)=\rm{dist}$$_g(x,p)$ for some fixed point $p$ in $\O$.

Similar as Reilly \cite{Reilly}, \eqref{xeq0} can be applied to reprove Alexandrov's theorem in $\mathbb{S}_+^n$ and $\mathbb{H}^n$, which is due to Alexandrov \cite{Al}. 
\begin{theorem}[Alexandrov, \cite{Al}]\label{Al}
Let $M$ be an embedded closed hypersurface in $\SS^n_+$ or $\HH^n$ with constant mean curvature $H$. Then $M$ must be a geodesic sphere.
\end{theorem}

Based on \eqref{xeq0}, we are also able to give an alternative proof of Theorem B. Moreover, our approach enables us to give a new Heintze-Karcher inequality for more general Riemannian manifolds with boundary.
\begin{theorem}\label{main theorem}
Let $(\Omega^{n}, g)$ be a $n$-dimensional compact  Riemannian manifold with smooth boundary $M$. Assume that the sectional curvature of $\O$ has a lower bound $\rm{Sect}\geq -1$. Let $H$ be the normalized mean curvature of $M$. Let $V(x)=\cosh r(x)$, where $r(x)=\rm{dist}$$_g(x,p)$ for some fixed point $p$ in $\O$. If $H$ is positive everywhere, then
 \begin{eqnarray}\label{ros1}
\int_M \frac{V}{H} dA  \geq\int_M \on_\nu V dA= \int_\O \ode V d \O.
\end{eqnarray}
The equality in \eqref{ros1} holds if and only if  $\O$ is a geodesic ball in a space form whose sectional curvature is $-1$.
\end{theorem}

Note that if $\Omega$ is of constant sectional curvature $-1$, then $\ode V =nV$. Then Theorem \ref{main theorem} reduces to Theorem B for the case $\O\subset\mathbb{H}^n$.
We remark that the proof of Theorem \ref{main theorem} also applies to the case $\O\subset \mathbb{S}_+^n$ to show \eqref{brendle}. Therefore, we give an alternative proof of Theorem B.

Inequality \eqref{ros1} is motivated by Brendle's inequality for warped product spaces \cite{Br}. However, except for the case of constant curvature manifolds, they are not the same.  First,  in Theorem \ref{main theorem}, we assume a lower bound of sectional curvature for Riemannian manifolds. Brendle \cite{Br} assumed the warped product structure for Riemannian manifolds and some conditions on the warped product function. Second, the equality case in our inequality can only occur for constant curvature manifolds. The equality case in Brendle's inequality can occur for the warped product spaces he considered. 

Comparing Theorem \ref{main theorem} with Theorem A, we want to ask whether Theorem \ref{main theorem} holds if  $\rm{Ric}$$\geq -(n-1)g$? Also if the right hand side of \eqref{ros1} is replaced by  $n \int_\O V d \O$, is \eqref{ros1} still true? Note that for $\rm{Ric}$$\geq -(n-1)g$, we have $\ode V \leq nV$. Hence the inequality
$\int_M \frac{V}{H} dA  \geq n\int_\O V d \O$
 is stronger than \eqref{ros1}.

Our method is based on the general Reilly's formula \eqref{xeq0} and the following Dirichlet
boundary value problem
\begin{equation}\label{op}
\left\{
\begin{array}{rccl}
\ode f&=&n f&
\hbox{ in } \Omega,\\
f&=&c  &\hbox{ on } M,\\
\end{array}
\right.
\end{equation}
for some real constant $c>0$. The existence and the regularity of solutions to \eqref{op} follows from standard theory of second order elliptic PDEs.

The paper is organized as follows. Section 2 is devoted to prove the general Reilly's formula \eqref{xeq0}. In Section 3, we use \eqref{xeq0}  to give an alternative proof of Theorem \ref{Al}. In Section 4, we use \eqref{xeq0} to prove the new Heintze-Karcher  inequality, Theorem \ref{main theorem}.


\section{General Reilly's formula}

For simplicity, We use $f_i, f_{ij},\cdots$ and $f_\nu$ to denote covariant derivatives and normal derivative of a function $f$ with respect to $g$ respectively. 

By integration by parts and Ricci identity, we have
\begin{eqnarray}\label{xeq1}
&&\int_\O V|\on^2 f|^2 d\O=\int_\O V\sum^{n}_{i,j=1}f_{ij}f_{ij} d\O \nonumber\\&=&\int_M V\sum^n_{i=1}f_{i\nu}f_i dA-\int_\O \sum^n_{i,j=1} V_jf_{ij}f_id\O-\int_\O V\sum^n_{i,j=1} f_{ijj}f_i  d\O\nonumber\\&=&\int_M V\sum^n_{i=1} f_{i\nu}f_i dA-\int_\O \sum^n_{j=1} V_j(\frac12|\on f|^2)_j d\O\nonumber\\&&-\int_\O \sum^n_{i=1} V\left((\ode f)_i+\sum^n_{j=1} R_{ij}f_j\right)f_i  d\O\nonumber
\\&=&\int_M V \sum^{n}_{i=1} f_{i\nu}f_i dA-\int_M \frac12|\on f|^2V_\nu  dA+\int_\O \frac12|\on f|^2\ode Vd\O\nonumber\\&&-\int_M V\ode f f_\nu dA+ \int_\O V(\ode f)^2 d\O+\int_\O \ode f \sum^n_{i=1} V_i f_i d\O\nonumber\\&&-\int_\O V \sum^n_{i,j=1}R_{ij}f_if_j d\O.
\end{eqnarray}
We also have
\begin{eqnarray}\label{xeq2}
\int_\O Vf\ode f d\O=\int_M Vff_{\nu} dA-\int_\O (V|\on f|^2+\sum^n_{i=1} V_if_i f) d\O.
\end{eqnarray}

Using \eqref{xeq1} and \eqref{xeq2}, we obtain
\begin{eqnarray}\label{xeq3}
&&\int_\O V((\ode f+Knf)^2-|\on^2 f+Kfg|^2)d\O\nonumber\\&=&\int_\O V((\ode f)^2-|\on^2 f|^2) d\O\nonumber\\&&+(2n-2)K\int_\O Vf\ode f d\O+n(n-1)K^2\int_\O Vf^2 d\O\nonumber\\&=&\int_M V\ode ff_\nu+\frac12|\on f|^2 V_\nu-V\sum^n_{i=1} f_{i\nu}f_i+(2n-2)KVff_{\nu} dA\nonumber\\&&+\int_\O -\frac12|\on f|^2\ode V - \ode f  \sum^n_{i=1} V_i f_i + V \sum^n_{i,j=1}R_{ij}f_if_j d\O\nonumber\\&&-(2n-2)K\int_\O (V|\on f|^2+\sum^n_{i=1} V_if_i f )d\O+n(n-1)K^2 \int_\O Vf^2 d\O.
\end{eqnarray}

We deal with the terms $\int_\O -\ode f\sum\limits^n_{i=1} V_i f_i  d\O$ and $-(2n-2)K\int_\O \sum\limits^n_{i=1} V_if_i f d\O$ in \eqref{xeq3} by integration by parts again.
\begin{eqnarray}\label{xeq4}
\int_\O -\ode  f\sum\limits^n_{i=1} V_i f_i  d\O&=&\int_M -f_\nu \sum\limits^n_{i=1} V_if_i dA+\int_\O \sum\limits^n_{i,j=1} V_{ij}f_if_j+ \sum\limits^n_{i=1} V_i(\frac12|\on f|^2)_i d\O\nonumber\\&=&\int_M -f_\nu \sum\limits^n_{i=1} V_if_i +\frac12|\on f|^2V_\nu dA\nonumber\\&&+\int_\O \sum\limits^n_{i,j=1} V_{ij}f_if_j-\frac12\ode V |\on f|^2 d\O.
\end{eqnarray}
\begin{eqnarray}\label{xeq5}
\int_\O \sum\limits^n_{i=1} V_if_i f d\O=\int_\O \sum\limits^n_{i=1} V_i(\frac12f^2)_id\O=
\int_M \frac12f^2V_\nu dA-\int_\O \frac12f^2\ode V d\O.
\end{eqnarray}

Inserting \eqref{xeq4} and \eqref{xeq5} into \eqref{xeq3}, we obtain
\begin{eqnarray}\label{xeq6}
&&\int_\O V((\ode f+Knf)^2-|\on^2 f+Kfg|^2)d\O\nonumber\\&=&\int_M V\ode ff_\nu+|\on f|^2 V_\nu-V\sum\limits^n_{i=1} f_{i\nu}f_i+(2n-2)KVff_{\nu}\nonumber\\&&\quad\quad-f_\nu \sum\limits^n_{i=1} V_if_i-(n-1)Kf^2V_\nu dA\nonumber\\&&+\int_\O\sum\limits^n_{i,j=1} V_{ij}f_if_j-\ode V|\on f|^2 -(2n-2)KV|\on f|^2+V\sum\limits^n_{i,j=1} R_{ij}f_if_jd\O\nonumber\\&&+(n-1)\int_\O (K\ode V+K^2 n V)f^2 d\O.
\end{eqnarray}

We now handle the boundary term in \eqref{xeq6}. We choose an orthonormal frame $\{e_i\}_{i=1}^n$ such that $e_n=\nu$ on $M$. Note that $z=f|_M$ and $u=f_\nu$.
From Gauss-Weingarten formula we deduce
\begin{eqnarray}\label{xeq7}
&&\int_M V\ode ff_\nu-V\sum_{i=1}^n f_{i\nu}f_i dA=\int_M V\sum_{i=1}^{n-1}f_{ii}f_\nu-V\sum_{i=1}^{n-1}f_{i\nu}f_i dA\nonumber\\&=&\int_M V\left(u\De z+(n-1)Hu^2-\<\n u, \n z\>+h(\n z, \n z)\right)dA.
\end{eqnarray}

On the other hand,
\begin{eqnarray}\label{xeq8}
\int_M |\on f|^2 V_\nu-\sum_{i=1}^n f_\nu V_if_i dA&=&\int_M |\n z|^2 V_\nu-u\< \n V, \n  z\> dA\nonumber\\&=&\int_M |\n z|^2 V_\nu+V\<\n u, \n z\>+Vu\De z dA.
\end{eqnarray}
It follows from \eqref{xeq7} and \eqref{xeq8} that
\begin{eqnarray}\label{xeq9}
&&\int_M V\ode ff_\nu+|\on f|^2 V_\nu-V\sum_{i=1}^nf_{i\nu}f_i+(2n-2)KVff_{\nu}\nonumber\\&&\quad\quad-f_\nu \sum_{i=1}^n V_if_i-(n-1)Kf^2V_\nu dA\nonumber\\&=&\int_M V\left(2u\De z+(n-1)Hu^2+h(\n z, \n z)+(2n-2)Kuz\right)dA\nonumber\\&&+\int_M V_\nu\left(|\n z|^2-(n-1)Kz^2\right) dA.
\end{eqnarray}

Combining \eqref{xeq6} and \eqref{xeq9}, we arrive at \eqref{xeq0}. \qed

\section{Alternative proof of Alexandrov's Theorem}

In this section, we prove Alexandrov's Theorem in $\SS_+^n$ and $\HH^n$ by using the general Reilly's formula \eqref{xeq0}. 

For the case $\O\subset \SS_+^n$, we choose $K=1$ and  $V(x)=\cos r(x)$. For the case $\O\subset \HH^n$, we choose $K=-1$ and  $V(x)=\cosh r(x)$. In any case, we have $\on^2 V=-KVg$. The general Reilly's formula \eqref{xeq0} reduces to
\begin{eqnarray}\label{xxeq0}
&&\int_\O V\left((\ode f+Knf)^2-|\on^2 f+Kfg|^2\right)d\O\nonumber\\&=&\int_M V\left(2u\De z+(n-1)Hu^2+h(\n z, \n z)+(2n-2)Kuz\right)dA\nonumber\\&&+\int_M \on_\nu V\left(|\n z|^2-(n-1)Kz^2\right) dA.\end{eqnarray}

Let $f: \O\to \RR$ be the solution of \begin{equation*}
\left\{
\begin{array}{rccl}
\ode f+Knf&=&1&
\hbox{ in } \Omega,\\
f&=&0 &\hbox{ on } M.\\
\end{array}
\right.
\end{equation*}

Then from \eqref{xxeq0} and Schwarz's inequality we obtain
\begin{eqnarray}\label{xxeq1}
\frac{n-1}{n}\int_\O V d\O\geq \int_M (n-1)Hu^2VdA.\end{eqnarray}
Since $H$ is a constant, we have from \eqref{xxeq1} and H\"older inequality that
\begin{eqnarray}\label{xxeq2}
\int_\O V d\O\geq nH\int_M u^2VdA\geq nH\frac{\left(\int_M uVdA\right)^2}{\int_M VdA}.\end{eqnarray}
By Green's formula and $\ode V=-KnV$, we have
\begin{eqnarray}\label{xxeq3}
\int_M uVdA=\int_\O \ode fV-f\ode Vd\O=\int_\O Vd\O.\end{eqnarray}
On the other hand, Minkowski formula in $\SS_+^n$ or $\HH^n$ tells that
\begin{eqnarray}\label{xxeq4}
\int_M VdA=\int_M HpdA=H\int_M pdA=nH\int_\O Vd\O.\end{eqnarray}
Here $p=\sin r\<\on r, \nu\>$ for $\O\subset \SS_+^n$ and $p=\sinh r\<\on r, \nu\>$ for $\O\subset \HH^n$.

Combining \eqref{xxeq1}--\eqref{xxeq4}, we find equality holds in \eqref{xxeq1}. Thus we must have  $|\on^2 f+Kfg|^2=\frac1n(\ode f+Knf)^2$. Taking into account that $\ode f+Knf=1$, we obtain $$\on^2 (f+\frac1n)=-K(f+\frac1n)g\hbox{ in }\O.$$
Since $f+\frac1n|_{M}=\frac1n$, it follows from an Obata type result (see Reilly \cite{Reilly2}) that $\O$ must be a geodesic ball. We complete the proof of Theroem \ref{Al}.

\section{Heintze-Karcher type  inequality}

In this section, we prove  Theorem \ref{main theorem} by using the general Reilly's formula \eqref{xeq0}. 

We let $K=-1$ and $V(x)=\cosh r(x)$ in  \eqref{xeq0},  $f$ be the solution to the following  Dirichlet boundary value problem:\begin{equation*}
\left\{
\begin{array}{rccl}
\ode f&=&n f&
\hbox{ in } \Omega,\\
f&=&c>0  &\hbox{ on } M.\\
\end{array}
\right.
\end{equation*}

By assumption, $\rm{Sect}$ $\geq -1$,  Hessian comparison theorem (see \cite{SY}) tells that  $r(x)$ $=\rm{dist}$$(x,p)$ satisfies \begin{eqnarray}\label{xxeq}
\nabla^2 r\leq \coth  r(g-dr^2),\quad \hbox{ for }x\in \O\setminus Cut(p),
\end{eqnarray}
where Cut$(p)$ is the cut locus of $p$. Thus
\begin{align}\label{xeq10}
  \nabla^2 V= \sinh r\n^2 r+\cosh r dr^2\leq \cosh r g=Vg, \quad \hbox{ for }x\in \O\setminus Cut(p). \end{align}
It follows that
\begin{align}\label{xeq11}
  \De V g -  \nabla^2 V \leq & (n-1) V g, \quad   \De V \leq n V, \quad \hbox{ for }x\in \O\setminus Cut(p).
\end{align}
Since Cut$(p)$ has zero measure, we obtain from \eqref{xeq11} and $Ric\geq -(n-1)g$ that
\begin{eqnarray}\label{xeq12}
\int_\O \left(\on^2 V-\ode Vg+(2n-2)Vg+V Ric\right)(\on f, \on f)d\O\geq 0,\end{eqnarray}
and
\begin{eqnarray}\label{xeq13}
(n-1)\int_\O (-\ode V+n V)f^2 d\O\geq 0.\end{eqnarray}
Using \eqref{xeq0}, \eqref{xeq12}, \eqref{xeq13}, Schwarz's inequality and $f|_M=c$, we have
\begin{eqnarray}\label{xeq14}
0&=&\frac{n-1}{n}\int_\O  V(\ode f-nf)^2 d\O\nonumber\\
&\geq & \int_\O V\left((\ode f-nf)^2-|\on^2 f-fg|^2\right)d\O\nonumber\\&\geq &\int_M (n-1)Hu^2V-(2n-2)cuV+(n-1)c^2\on_\nu VdA.
\end{eqnarray}

By H\"older inequality and \eqref{xeq14} we deduce that
\begin{align}\label{xeq15}
  \left(\int_{M} u V dA\right)^2 \leq & \int_{M} \frac{V}{H} dA \int_{M} H u^2V dA \nonumber\\
  \leq & \int_{M} \frac{V}{H} dA \int_{M} 2 cuV   -  c^2 \on_\nu V dA.
\end{align}
It follows that
\begin{align}\label{xeq16}
  &\left(\int_{M} uV dA - c \int_{M} \frac{V}{H} dA\right)^2 - c^2 \left(\int_{M} \frac{V}{H} dA\right)^2 \nonumber \\
   & \leq -c^2 \int_{M} \frac{V}{H} dA \int_{M} \on_\nu V  dA.
\end{align}
Thus
\begin{align}\label{xeq17}
- c^2 \left(\int_{M} \frac{V}{H} dA\right)^2 \leq -c^2 \int_{M} \frac{V}{H} dA \int_{M} \on_\nu V  dA.
\end{align}
Since $c\neq 0$, by eliminating $-c^2\int_{M} \frac{V}{H} dA$ from both sides of \eqref{xeq14}, we conclude
\begin{align}\label{xeq18}
  \int_{M} \frac{V}{H} dA & \geq \int_{M} \on_\nu V  dA = \int_\O \ode Vd\O.
\end{align}

We now explore the equality case. If the equality holds in \eqref{xeq18}, then equality must hold in \eqref{xeq14}. It follows that $|\on^2 f-fg|^2=\frac1n(\ode f-nf)^2$. Taking into account of $\ode f=nf$, we obtain $$\on^2 f=fg\hbox{ in }\O.$$
Since $f|_{M}=c$, it follows from an Obata type result (see Reilly \cite{Reilly2}) that $\O$ must be a geodesic ball in a space form with constant sectional curvature $-1$.
We complete the proof of Theorem \ref{main theorem}.

\

For the case $\O\subset\mathbb{S}_+^n$, our method still applies to prove \eqref{brendle}, that is
\begin{eqnarray}\label{xeq19}
\int_M \frac{\cos r}{H}dA\geq n\int_\Omega \cos rd\O.
\end{eqnarray}
Since the proof is almost the same, we only indicate the difference.
First, we choose $K=1$ and $V(x)=\cos r(x)>0$ in \eqref{xeq0} and $f$ be the solution to the Dirichlet boundary value problem:\begin{equation*}
\left\{
\begin{array}{rccl}
\ode f&=&-n f&
\hbox{ in } \Omega,\\
f&=&c>0  &\hbox{ on } M.\\
\end{array}
\right.
\end{equation*}
We notice that $\on^2V=-Vg$ and $\ode V=-nV$. Then the same argument as above works to show \eqref{xeq16}. \qed

\

\noindent\textbf{Acknowledgments.}  The research of the first author is supported by National Natural Foundation of China (Grant No. 11301497). The second author has received funding from the European Research Council under the European Union's Seventh Framework Programme (FP7/2007-2013) / ERC grant
agreement no. 267087.  Both authors would like to express their gratitude to Prof. Xinan Ma and Prof. Guofang Wang for their constant encouragement and support.

\

\end{document}